\documentclass[11pt, reqno]{amsart}

\usepackage[left=3.5cm, right=3.5cm, top=2.5cm, bottom=2.5cm]{geometry}

\usepackage[utf8]{inputenc}

\usepackage{amssymb, amsmath, amsthm, esint}

\usepackage{latexsym}

\usepackage{graphicx}

\usepackage{color}

\newtheorem{defi}{Definition}[section]
\newtheorem{thm}[defi]{Theorem}
\newtheorem{ex}[defi]{Example}
\newtheorem{rem}[defi]{Remark}
\newtheorem{prop}[defi]{Proposition}
\newtheorem{lemma}[defi]{Lemma}
\newtheorem{cor}[defi]{Corollary}

\newcommand{\Sp}{\mathbb S}
\newcommand{\R}{\mathbb R}
\newcommand{\N}{\mathbb N}
\newcommand{\V}{{\rm Vol}}

\title[Steklov eigenvalues of submanifolds]{Upper bounds for Steklov eigenvalues of
 submanifolds in Euclidean space via the intersection index}
\author{Bruno Colbois}
\address{Universit\'e de Neuch\^atel, Institut de Math\'ematiques, Rue
  Emile-Argand 11, CH-2000 Neuch\^atel, Switzerland}
\email{bruno.colbois@unine.ch}
\author{Katie Gittins}
\address{Department of Mathematical Sciences, Durham University, Lower Mountjoy, DH1 3LE Durham,
United Kingdom.}
\email{katie.gittins@durham.ac.uk}
\subjclass[2010]{35P15, 58C40}
\keywords{Steklov problem, submanifold, Euclidean space, intersection index}
\date{\today}

\begin{document}
\begin{abstract}
We obtain upper bounds for the Steklov eigenvalues $\sigma_k(M)$ of a smooth, compact, connected,
$n$-dimensional submanifold $M$ of Euclidean space with boundary $\Sigma$ that involve the intersection
indices of $M$ and of $\Sigma$. One of our main results is an explicit upper bound in terms of the intersection index of $\Sigma$, the volume of $\Sigma$ and the volume of $M$ as well as dimensional constants.
By also taking the injectivity radius of $\Sigma$ into account, we obtain an upper bound that has the optimal exponent of $k$ with respect to the asymptotics of the Steklov eigenvalues as $k \to \infty$.
\end{abstract}
\maketitle

\vspace{-0.5cm}
\section{Introduction}
Let $(M,g)$ be a smooth, compact, connected
Riemannian manifold of dimension $n\geq 2$ with smooth boundary $\Sigma$. The Steklov eigenvalues of $M$  are the real numbers $\sigma$ for which there exists a non-zero harmonic function $f:M\rightarrow\R$ that satisfies $\partial_\nu f=\sigma f$ on $\Sigma$. Here and in what follows, $\partial_\nu$ denotes the outward-pointing normal derivative on $\Sigma$.
It is well known that the Steklov eigenvalues can be written in a non-decreasing sequence
$$0=\sigma_0(M)<\sigma_1(M)\leq\sigma_2(M)\leq\cdots,$$
where each eigenvalue is repeated according to its multiplicity, and the only point of accumulation is $+\infty$.
The Steklov eigenvalues satisfy the following asymptotic formula
\begin{equation}\label{eq:1.1}
  \sigma_k(M) = 2\pi \left( \frac{k}{\omega_{n-1} \vert \Sigma \vert} \right)^{1/(n-1)} + O(1), \quad k \to \infty,
\end{equation}
where $\omega_{n-1}$ is the volume of a ball of radius 1 in $\R^{n-1}$, see for example \cite{GP14} and references therein.

In recent years, the interplay between the Steklov eigenvalues of a manifold and the geometry of
the manifold and its boundary has been intensively studied in the field of spectral geometry.
One way to shed light on this interplay is to obtain upper bounds for the Steklov eigenvalues in terms of some of the geometric quantities of $M$ and $\Sigma$ .

An interesting and natural problem is to investigate in which geometric situations it is possible to obtain upper bounds for the Steklov eigenvalues where the exponent of $k$ is optimal with respect to the Steklov spectral asymptotics \eqref{eq:1.1}.

A compelling result in this direction is due to Provenzano and Stubbe \cite{PS} which asserts that
for a bounded open set $\Omega \subset \R^{n+1}$ with connected boundary $\partial \Omega$ of class $C^2$,
for $k \in \N$,
\begin{gather*}
  \lambda_k(\partial \Omega) \leq \sigma_k(\Omega)^2 + 2 c_{\Omega} \sigma_k(\Omega), \\
  \sigma_k(\Omega) \leq c_{\Omega} + \sqrt{c_{\Omega}^2 + \lambda_k(\partial \Omega)}, \\
  \vert \sigma_k(\Omega) - \lambda_k(\partial \Omega) \vert \leq c_{\Omega},
\end{gather*}
where $\lambda_k(\partial \Omega)$ denotes the $k$-th eigenvalue of the Laplacian on $\partial \Omega$
and $c_{\Omega}$ is a constant that depends on the dimension, the maximum of the mean of the absolute values
of the principal curvatures on $\partial \Omega$ and the rolling radius of $\Omega$. We observe that due to
the Weyl Law, the exponent of $k$ in these bounds agrees with the asymptotics of the Steklov eigenvalues \eqref{eq:1.1}.

These results were generalised to a Riemannian setting of an open set with convex boundary in a complete
Riemannian manifold with either non-positive or strictly positive sectional curvature and with some
hypotheses on the second fundamental form by Xiong \cite{X}. Analogous results were obtained in \cite{CGH} in
the more general setting of smooth, compact Riemannian manifolds with boundary under assumptions
on the sectional curvatures in a tubular neighbourhood of the boundary, the principal curvatures of the boundary,
and the rolling radius.

In \cite{CGG}, we also obtained bounds for the Steklov eigenvalues that have the optimal exponent of $k$
in the special case of hypersurfaces of revolution in Euclidean space (see Theorem 1.8, Theorem 1.11 and Proposition 3.1).

The first goal of this article is to obtain an upper bound for the Steklov eigenvalues of a submanifold of $\R^m$ with non-empty boundary which has optimal exponent of $k$ with respect to \eqref{eq:1.1} and does not depend on the curvatures of $M$ or of $\Sigma$. This is achieved in Theorem \ref{thm:3} below.

On the other hand, we note that there are also important results where the exponent of $k$ is not optimal.
In the Riemannian setting, it was shown by Colbois, El Soufi and Girouard \cite{ceg1} that if $M$ is an $n$-dimensional Riemannian manifold that is conformally equivalent to a complete Riemannian manifold with non-negative Ricci curvature, then for any domain $\Omega \subset M$ and $k \in \N$,
\begin{equation} \label{eq:CEG}
  \sigma_k(\Omega) \leq \alpha(n) \frac{\vert\Omega\vert_g^{(n-2)/n}}{\vert \partial \Omega\vert} k^{2/n},
\end{equation}
where $\alpha(n) >0$ is a constant that depends only on $n$. It is also possible to express this inequality in terms of the isoperimetric ratio $I(\Omega)=\frac{\vert\partial \Omega \vert}{\vert \Omega\vert_g^{(n-1)/n}}$ of $\Omega$
as follows:
\begin{equation} \label{iso}
  \sigma_k(\Omega) \vert\partial \Omega\vert^{1/(n-1)} \leq \alpha(n) \frac{ k^{2/n}}{I(\Omega)^{(n-2)/(n-1)}}  .
\end{equation}
In addition \cite{AH11} addresses the case of a Riemannian manifold where the Ricci curvature is bounded from below by a negative constant.

For a general submanifold $M$ in Euclidean space with fixed boundary $\Sigma$, we obtained the following upper bound for the Steklov eigenvalues (see Theorem 1.1 of \cite{CGG}):

\begin{thm}\label{thm:CGG}
  Let $n \geq 2$. Let $\Sigma$ be a fixed $(n-1)$-dimensional compact, smooth submanifold in $\R^{m}$.
  There exists a constant $A_\Sigma$ depending on $\Sigma$ such that any compact, $n$-dimensional
  submanifold $M$ of $\R^m$ with boundary $\Sigma$ satisfies
 \begin{equation}\label{eq:CGG}
   \sigma_k(M) \le A_\Sigma \vert M\vert k^{2/(n-1)}.
 \end{equation}
\end{thm}

The constant $A_{\Sigma}$ depends upon several geometric quantities of $\Sigma$ including the Ricci curvature, the diameter and the distortion, and some of these quantities are difficult to measure.

In fact, in Section \ref{s:5}, we observe that, in general, it is not true that a bound of the form
\begin{equation*}
\sigma_k(M) \leq C(\vert \Sigma \vert, n)\vert M \vert^{\beta} k^{\alpha}
\end{equation*}
with $\beta >0$ has $\alpha = \frac{1}{n-1}$ (where $C(\vert \Sigma \vert, n)$ is a constant depending on $\vert \Sigma \vert$, $n$). Indeed, for a cylinder $M = \Sigma \times [0,L]$, if $\beta >0$ then $\alpha > \frac{1}{n-1}$.

The second goal of this article is to obtain estimates of the same form as \eqref{eq:CGG} and \eqref{iso} that depend on an invariant that is more robust than the Ricci curvature and more easily computed and understood than the distortion. This is addressed in Theorem \ref{thm:1}.

The invariant that we consider throughout this paper is the intersection index which we define below as in \cite{CDES}.
For a compact immersed submanifold $N$ of dimension $q$ in $\R^{q+p}$, almost every
$p$-plane $\Pi$ in $\R^{q+p}$ is transverse to $N$ so that the intersection $\Pi \cap N$
consists of a finite number of points.
\begin{defi}
\label{def:index}
The intersection index of $N$ is
\begin{equation}
\label{eq:index}
i_p(N) := \sup_{\Pi} \{\# \Pi \cap N\},
\end{equation}
where the supremum is taken over the set of all $p$-planes $\Pi$ that are transverse to $N$ in $\R^{q+p}$.
\end{defi}
Our main result is the following.
\begin{thm}
\label{thm:3}
Let $n \geq 2$.
Let $\Sigma$ be an $(n-1)$-dimensional,
closed, smooth submanifold of $\R^m$. Let $r_0$ denote the injectivity radius of $\Sigma$.
Let $M$ be a compact, $n$-dimensional submanifold of $\R^m$ with boundary $\Sigma$.
There exist constants $\tilde{A}(n), \tilde{B}(n)>0$ depending only on $n$, such that for $k \geq 1$,
\begin{equation}\label{eq:T3}
  \sigma_k(M) \leq \tilde{A}(n) \frac{i(M)}{r_0} + \tilde{B}(n) \, i(M)\left(\frac{i(\Sigma) k}{\vert \Sigma \vert}\right)^{1/(n-1)},
\end{equation}
where $i(\Sigma) = i_{m-(n-1)}(\Sigma)$ and $i(M) = i_{m-n}(M)$.
\end{thm}

Upper bounds for the eigenvalues of the Laplace--Beltrami operator on a closed, connected submanifold of $\R^m$ in terms of the intersection index that are optimal with respect to the Weyl Law were obtained in \cite{CDES}.

\begin{rem}
In the upper bound (\ref{eq:T3}), there is a strong dependence on the injectivity radius $r_0$ of $\Sigma$. First, the term involving $r_0$ is separate from the term involving $k$, so that, asymptotically, we have the optimal behaviour in terms of the exponent of $k$. Moreover, it turns out that we need to take the injectivity radius $r_0$ into account. In Section \ref{s:4}, for $n \geq 3$, we construct an example of a compact $(n+1)$-dimensional submanifold $M$ of $\R^{n+3}$ which has $\vert \partial M \vert = 1$ and $i(\partial M)$, $i(M)$ bounded, but $\sigma_1(M) \to +\infty$ as the injectivity radius of $\partial M$ tends to zero. This shows that the injectivity radius is necessary in \eqref{eq:T3} when the dimension of $M$ is at least 4.
\end{rem}

For the special case where $\Sigma$ is a hypersurface, we obtain the following corollary which does not depend on $i(M)$.
\begin{cor}\label{cor:hyp}
Let $n \geq 2$.
Let $\Sigma \subset \R^n$ be an $(n-1)$-dimensional, closed, smooth hypersurface bounding a domain $M \subset \R^n$. Let $r_0$ denote the injectivity radius of $\Sigma$.
There exist constants $\hat{A}(n), \hat{B}(n) >0$ depending only on $n$ such that for $k \geq 1$,
\begin{equation*}
\sigma_k(M) \leq \frac{\hat{A}(n)}{r_0} + \hat{B}(n) \left(\frac{i(\Sigma) k}{\vert \Sigma \vert}\right)^{1/(n-1)}.
\end{equation*}
\end{cor}

This result is in the same spirit as Corollary 4.8 of \cite{PS} except that the constant term in the above bound depends only on the dimension and the injectivity radius of $\Sigma$ while the coefficient of
$k^{1/(n-1)}$ also depends on $i(\Sigma)$.
On the other hand, if we assume that $\Sigma$ is convex then $i(\Sigma) = 2$ so we separate the geometry from the asymptotics similarly to Corollary 5.4 of \cite{PS}.

By making use of the intersection index of $\Sigma$, we obtain the following more explicit version of Theorem \ref{thm:CGG} that does not depend on the Ricci curvature or the distortion of $\Sigma$.
\begin{thm}
\label{thm:1}
Let $n \geq 2$.
Let $\Sigma$ be an $(n-1)$-dimensional,
compact, smooth submanifold of $\R^m$.
Let $M$ be a closed, $n$-dimensional submanifold of $\R^m$ with boundary $\Sigma$.
There exists a constant $\tilde{C}(n,m)>0$ depending only on $n$, $m$ such that for $k \geq 1$,
\begin{equation}\label{eq:T1}
  \sigma_k(M) \leq \tilde{C}(n,m) \frac{i(\Sigma)^{2/(n-1)} \vert M \vert}{\vert \Sigma \vert^{(n+1)/(n-1)}} k^{2/(n-1)}.
\end{equation}
\end{thm}
This theorem gives rise to an upper bound for $A_{\Sigma}$ in terms of the intersection index and the volume of $\Sigma$ as well as a dimensional constant.

Another way to express Inequality \eqref{eq:T1} in Theorem \ref{thm:1} is as follows.
\begin{cor}
Under the same assumptions as in Theorem \ref{thm:1},
\begin{equation}\label{eq:}
  \sigma_k(M) \vert \Sigma \vert^{1/(n-1)} \leq \tilde{C}(n,m) \frac{i(\Sigma)^{2/(n-1)}  }{I(M)^{n/(n-1)}} k^{2/(n-1)},
\end{equation}
where $I(M)$ denotes the isoperimetric ratio $I(M)=\frac{\vert \Sigma \vert}{ \vert M\vert^{(n-1)/n}}$.
\end{cor}

Inequality \eqref{eq:} looks like Inequality \eqref{iso} except that the exponent $\frac{2}{n-1}$ of $k$ in the former is slightly worse than the exponent $\frac{2}{n}$ of $k$ in the latter. However, Inequality \eqref{iso} was established for domains of a complete manifold with Ricci curvature conformally non-negative, and Inequality \eqref{eq:} holds for any submanifold $\Sigma$ of Euclidean space without any curvature assumptions.

\begin{rem} Estimates as in Theorem \ref{thm:3} and Theorem \ref{thm:1}, in terms of the intersection indices of $M$ and $\Sigma$, are particularly interesting if we consider algebraic submanifolds $M$ and $\Sigma$ of $\mathbb R^m$. In this case, the intersection indices of $M$ and $\Sigma$ are bounded from above by the product of the degrees of the polynomials that define these submanifolds (see the proof of Corollary 4.1 of \cite{CDES} and references therein). We note that, in general, for a submanifold that is defined as the level set of a function it is not possible to estimate the curvature and the distortion.
\end{rem}

\subsection*{Organisation of the paper}
In Section \ref{s:2}, we recall some general results about the Steklov eigenvalue problem.
Section \ref{s:3} is devoted to the proofs of Theorem \ref{thm:3} and Theorem \ref{thm:1}.
In Section \ref{s:4}, we construct an example to show that the term involving the injectivity radius in Theorem \ref{thm:3} is necessary in the case where the dimension is at least 4. In Section \ref{s:5}, we present two
examples for which, if the upper bound for $\sigma_k(M)$ contains the volume of $M$ to a positive power, then the exponent of $k$ cannot be optimal with respect to the Steklov spectral asymptotics.

\section{Some general facts about the Steklov eigenvalue problem} \label{s:2}

For $k \in \N^*$, the Steklov eigenvalues of $(M,g)$ can be characterised by the following variational formula:
\begin{equation}
  \label{minmax}
  \sigma_k(M) = \min_{E \in \mathcal{H}_k} \max_{0 \neq u \in E} R_M(u),
\end{equation}
where $\mathcal{H}_k$ is the set of all $k$-dimensional subspaces in the Sobolev space $H^1(M)$
that are orthogonal to constants on $\Sigma$, and
\begin{equation*}
  \label{Rayleighquotient}
  R_M(u) = \frac{\int_{M} \vert \nabla u \vert^2 \, d\V_M}{\int_{\Sigma} \vert u \vert^2 \, d \V_{\Sigma}}
\end{equation*}
is the corresponding Rayleigh quotient.

For the first non-zero Steklov eigenvalue of $M$, we have that
\begin{equation*}
  \sigma_1(M) = \min \left\{ R_M(u): 0 \neq u \in H^{1}(M), \int_{\Sigma} u \, d\V_{\Sigma} = 0 \right\}.
\end{equation*}
\indent
In Section \ref{s:4}, we will make use of the comparison with the eigenvalues of the mixed Steklov--Neumann problem. Let $\Sigma$ be an $(n-1)$-dimensional, closed, smooth submanifold of $\R^m$. Let $M$ be an $n$-dimensional, compact submanifold of $\R^m$ with boundary $\Sigma$. We consider $\Omega \subset M$ such that $\Sigma \subset \partial \Omega$. We denote the intersection of $\partial \Omega$ with the interior of $M$ by $\partial_{I} \Omega$ and suppose that it is smooth.
As in \cite{CGG}, the mixed Steklov-Neumann problem on $\Omega$ is given by
\begin{gather*}
  \Delta u=0\mbox{ in } \Omega,\\
  \partial_\nu u=\sigma u\mbox{ on } \Sigma,\\
  \partial_\nu u=0 \mbox{ on }\ \partial_I \Omega.
\end{gather*}
The corresponding eigenvalues of this problem form a discrete sequence
$$0=\sigma_0^N(\Omega) \leq \sigma_1^N(\Omega) \leq \sigma_2^N(\Omega)\leq\cdots,$$
where the only point of accumulation is $+\infty$.
For each $k \in \N^*$ the $k$-th eigenvalue is given by
\begin{equation*}
 \sigma_k^{N}(\Omega)=\min_{E\in\mathcal{H}_{k}(\Omega)}\max_{0\neq u\in E}\frac{\int_{\Omega}|\nabla u|^2\,d\V_M}{\int_{\Sigma}|u|^2\,d\V_\Sigma},
\end{equation*}
where $\mathcal{H}_{k}(\Omega)$ is the set of all $k$-dimensional subspaces in the Sobolev space $H^1(\Omega)$ which are orthogonal to constants on $\Sigma$. We therefore have that for $k \in \N^*$,
\begin{equation}\label{eq:SN}
  \sigma_k^{N}(\Omega) \leq \sigma_k(M).
\end{equation}

\section{Upper bounds for submanifolds in Euclidean space}
\label{s:3}

Let $n \geq 2$. Let $\Sigma$ be an $(n-1)$-dimensional, closed, smooth submanifold of $\R^m$.
Let $M$ be a compact, $n$-dimensional submanifold of $\R^m$ with boundary $\Sigma$.
In order to obtain upper bounds for the Steklov eigenvalues of $M$, we define test functions
with disjoint support and then estimate the Rayleigh quotient (see also \cite{CGG} and references therein).

We recall the following proposition from \cite{CDES} for a compact, immersed submanifold $N$ of dimension $q$ in $\R^{q+p}$. It shows that if the intersection index of $N$ is bounded then $N$ does not concentrate in small regions of $\R^{q+p}$.
\begin{prop}
\label{prop:noconcentration}
For all $x \in \R^{q+p}$ and all $r > 0$, we have that
\begin{equation}\label{eq:noconcentration}
  \vert  N \cap B(x,r) \vert \leq \frac{i_p(N)}{2} \vert \Sp^{q} \vert \, r^{q}.
\end{equation}
\end{prop}
We observe that if $\Sigma$ is an $(n-1)$-dimensional, compact, smooth submanifold of $\R^{m}$ such that $i_{m-(n-1)}(\Sigma)$ is bounded, then
\begin{equation}\label{eq:sigma}
 \vert \Sigma \cap B(x,r)\vert \leq \frac{i_{m-(n-1)}(\Sigma)}{2} \vert \Sp^{n-1} \vert \, r^{n-1}.
\end{equation}
In addition, if $M$ is an $n$-dimensional, compact, smooth submanifold of $\R^{m}$ such that $i_{m-n}(M)$ is bounded, then
\begin{equation}\label{eq:m}
 \vert M \cap B(x,r)\vert \leq \frac{i_{m- n}(M)}{2} \vert \Sp^{n} \vert \, r^{n}.
\end{equation}
We make use of Inequality \eqref{eq:sigma} in the proofs of Theorems \ref{thm:1} and \ref{thm:3}.
On the other hand, by also taking the injectivity radius into account and making use of Inequality \eqref{eq:m},
we obtain the upper bound with optimal exponent given in Theorem \ref{thm:3}.

\subsection{Proof of Theorem \ref{thm:1}}
\label{section:proofT1}
To obtain disjoint sets on which to define test functions, we apply a variation of Lemma 2.2 and Corollary 2.3 of \cite{cm} which appeared in \cite{CGG} (Lemma 4.1). We recall this result below.

\begin{lemma}\label{lem:CMrevisited}
  Let $(X,d,\mu)$ be a complete, locally compact metric measured space,
  where $\mu$ is a non-atomic finite measure. Assume that for all
  $r>0$, there exists an integer $C>0$ such that each ball of radius
  $r$ can be covered by $C$ balls of radius $r/2$.
  Let $K>0$. If there exists a radius $r>0$ such that, for each $x \in X$
  $$
  \mu(B(x,r)) \leq \frac{\mu(X)}{4C^2K},
  $$
  then, there exist $\mu$-measurable subsets $A_1,\dots,A_K$ of $X$
  such that, for all $i\le K$, $\mu(A_i)\ge \frac{\mu(X)}{2CK}$
  and, for
  $i\not =j$, $d(A_i,A_j) \geq 3r$.
\end{lemma}

\begin{proof}[Proof of Theorem \ref{thm:1}]
We note that in the case under consideration, the ambient space is $X = \R^m$ and we can choose
$C(m)=32^m$ (see for example \cite{ceg0, ColboisMTL, CGG} for further details).

The measure $\mu$ is defined for a Borelian set $\mathcal{O}$ of $\R^m$ as
$$\mu(\mathcal O)=\int_{\Sigma \cap \mathcal O} d \V_\Sigma =\vert \Sigma \cap \mathcal O\vert ,$$
and $\mu (\Sigma)$ is the usual volume $\vert \Sigma \vert$ of $\Sigma$ (see also \cite{ceg0, CGG}).

As we wish to begin with $2k+2$ test functions in what follows, we take $K = 2k+2$.

We choose
\begin{equation}\label{choice}
r = \left(\frac{\vert \Sigma \vert}{2C(m)^2 \vert \Sp^{n-1} \vert i(\Sigma) (2k+2)}\right)^{1/(n-1)}.
\end{equation}
By Proposition \ref{prop:noconcentration}, we have that for all $x \in \R^m = \R^{n-1 + p}$,
\begin{equation*}
  \mu(B(x,r)) = \int_{\Sigma \cap B(x,r)} \, d \V_{\Sigma}
  = \vert \Sigma \cap B(x,r)\vert \leq \frac{i(\Sigma)}{2} \vert \Sp^{n-1} \vert \, r^{n-1}.
\end{equation*}
Together with the above choice of $r$, we obtain
\begin{equation*}
  \mu(B(x,r)) \leq \frac{\vert \Sigma \vert}{4C(m)^2 (2k+2)}.
\end{equation*}
Hence by Lemma \ref{lem:CMrevisited}, there exist $2k+2$ $\mu$-measurable subsets $A_1, \dots, A_{2k+2}$ of $\R^m$ such that for all $i=1, \dots, 2k+2$,
\begin{enumerate}
\item[(i)] $\mu(A_i)\geq \frac{\vert \Sigma \vert}{2C(m)(2k+2)}$,
\item[(ii)] for $i \neq j$, $d(A_i,A_j) \geq 3r$.
\end{enumerate}

For each $i \in \{1, \dots, 2k+2\}$, we define the $r$-neighbourhood of $A_i$ as
$$ A_i^{r} := \{x \in \R^m : d(x,A_i) < r\},$$
where $d$ denotes the Euclidean distance.
By property (ii) above, the $A_i^{r}$ are disjoint.

Similarly to \cite{cm}, for each $i \in \{1, \dots, 2k+2\}$, we construct a test function $g_i$ with support in $A_i^{r}$ as follows. For $x \in A_i^{r}$,
\begin{equation*}
  g_i(x) = 1 - \frac{d(x,A_i)}{r}.
\end{equation*}

Then $\vert \nabla g_i \vert^2 \leq \frac{1}{r^2}$ almost everywhere in $A_i^r$. So
\begin{equation*}
  \int_{M \cap A_i^r} \vert \nabla g_i \vert^2 \, d \V_{M} \leq \frac{\vert M \cap A_i^{r} \vert}{r^2}.
\end{equation*}

Since $g_i(x) = 1$ for $x \in A_i$, we have that
\begin{equation*}
  \int_{\Sigma} g_i^2 \, d\V_{\Sigma} \geq \int_{\Sigma \cap A_i} d \V_{\Sigma}
  = \mu(A_i) \geq \frac{\vert \Sigma \vert}{2C(m)(2k+2)}.
\end{equation*}

We now use the previous two inequalities to estimate the Rayleigh quotient and hence prove Theorem \ref{thm:1}.

As the $A_i^r$ are disjoint, there exist $k+1$ of them, say $A_1^r, \dots, A_{k+1}^r$, such that for $i=1, \dots, k+1$,
\begin{equation*}
  \vert M \cap A_i^r \vert \leq \frac{\vert M \vert}{k+1}.
\end{equation*}
Therefore, using \eqref{choice}, we have that
\begin{align*}
  \frac{\int_{M} \vert \nabla g_i \vert^2 \, d \V_{M}}{\int_{\Sigma} g_i^2 \, d \V_{\Sigma}}
  &\leq  \frac{\vert M \vert}{(k+1)r^2} \frac{2C(m)(2k+2)}{\vert \Sigma \vert} \\
  &= 4 \cdot 4^{2/(n-1)} C(m)^{(n+3)/(n-1)} \vert \Sp^{n-1}\vert^{2/(n-1)} \frac{i(\Sigma)^{2/(n-1)} \vert M \vert}{\vert \Sigma \vert^{(n+1)/(n-1)}} (k+1)^{2/(n-1)}\\
  &\leq 4 \cdot 4^{3/(n-1)} C(m)^{(n+3)/(n-1)} \vert \Sp^{n-1} \vert^{2/(n-1)} \frac{i(\Sigma)^{2/(n-1)} \vert M \vert}{\vert \Sigma \vert^{(n+1)/(n-1)}} k^{2/(n-1)},
\end{align*}
which proves Theorem \ref{thm:1} with $\tilde{C}(n,m) = 4 \cdot 4^{3/(n-1)} C(m)^{(n+3)/(n-1)} \vert \Sp^{n-1} \vert^{2/(n-1)}$.
\end{proof}

\subsection{Proof of Theorem \ref{thm:3}}
\label{section:proofT3}
The intuition behind the proof of this theorem is that we can make use of Inequality \eqref{eq:m} instead of taking the whole volume of $M$ into account.
Indeed, Inequality \eqref{eq:m} gives local control on the volume of $M$ via the intersection index $i(M)$.
This allows us to obtain that the numerator in the Rayleigh quotient is bounded from above
by a geometric constant times $r^{n-2}$. If we also have that the denominator is bounded from below by a geometric constant times $r^{n-1}$, then this would give rise to an upper bound for the Rayleigh quotient as a geometric constant times the desired factor $\frac{1}{r}$ as opposed to $\frac{1}{r^2}$. By taking the injectivity radius into account, we can obtain a suitable lower bound for the denominator of the Rayleigh quotient.

In fact, in the case where we take the injectivity radius into account, we can obtain a collection of disjoint
balls on which to define the test functions. That is, we do not need to appeal to the result of \cite{cm}.
To do this we make use of inequalities \eqref{eq:sigma} and \eqref{eq:m} and the following result due
to Croke \cite{Croke}.

For the submanifold $M$ in $\R^m$, we denote by $g$ the Riemannian metric induced by $\R^m$ on $M$. Then, we have that the Euclidean distance between two points $d$ is smaller than or equal to the Riemannian distance $d_{g}$ between these points. So, for $r>0$, $x \in \Sigma$,
\begin{equation*}
 B_g(x,r) =\{y \in \Sigma : d_g(x,y) \leq r\} \subset B(x,r)=\{y \in \Sigma : d(x,y) \leq r\},
\end{equation*}
hence $\vert B(x,r) \vert \geq \vert B_g(x,r) \vert$.
By Proposition 14 of \cite{Croke}, for $r \leq \frac{r_0}{2}$, we have that there exists a constant $D_n > 0$, depending only on $n$, such that
\begin{equation}\label{eq:croke}
  \vert B(x_j,r) \cap \Sigma \vert \geq \vert B_g(x_j,r) \cap \Sigma \vert \geq \frac{2^{n-1} \vert \Sp^{n-1} \vert^n}{(n-1)^{n-1} \vert \Sp^n \vert^{n-1}} \, r^{n-1}.
\end{equation}

\begin{proof}[Proof of Theorem \ref{thm:3}]
We choose
\begin{equation}\label{eq:rchoice}
r = \left(\frac{\vert \Sigma \vert}{4^{n-1} \vert \Sp^{n-1} \vert \, i(\Sigma) \, k}\right)^{1/(n-1)}.
\end{equation}
Note that this choice of $r$ does not depend on $m$, the dimension of the ambient space, in contrast to the choice in \eqref{choice}.

Let $\{x_j\}_{j=1}^{N}$ be a maximal subset of points in $\Sigma$ such that $d(x_i,x_j) \geq 4r$ for $i \neq j$.
By maximality, we have that $\Sigma \subset \cup_{i=1}^{N} B(x_i, 4r)$.
Hence by \eqref{eq:sigma} we have
\begin{equation*}
  \frac{i(\Sigma)}{2} \vert \Sp^{n-1} \vert \, N (4r)^{n-1}
  \geq \sum_{i=1}^{N} \vert B(x_j, 4r) \cap \Sigma \vert \geq \vert \Sigma \vert,
\end{equation*}
which implies that
\begin{equation}\label{eq:lbdN}
  N \geq \frac{2 \vert \Sigma \vert}{4^{n-1} \, i(\Sigma) \, \vert \Sp^{n-1} \vert \, r^{n-1}}.
\end{equation}
To estimate $\sigma_k$, we require $k+1$ disjoint open balls.
In what follows, we will take the balls $B(x_i, 2r)$ which are disjoint since $d(x_i,x_j) \geq 4r$, $i \neq j$.
Thus we need $N \geq k+1$ so, by \eqref{eq:lbdN}, it suffices to take
\begin{equation*}
  \frac{2 \vert \Sigma \vert}{4^{n-1} \, i(\Sigma) \, \vert \Sp^{n-1} \vert \, r^{n-1}} >k.
\end{equation*}
By the choice of $r$ in \eqref{eq:rchoice}, we have
\begin{equation*}
  \frac{2 \vert \Sigma \vert}{4^{n-1} \, i(\Sigma) \, \vert \Sp^{n-1} \vert}
  \, \frac{4^{n-1} \vert \Sp^{n-1} \vert \, i(\Sigma) \, k}{\vert \Sigma \vert} = 2k
\end{equation*}
which is indeed larger than $k$.

We set
\begin{equation}\label{eq:k0}
  k_0 = \left\lceil \frac{\vert \Sigma \vert}{2^{n-1} \, i(\Sigma) \, \vert \Sp^{n-1}\vert \, r_0^{n-1}} \right\rceil
\end{equation}
and consider the case $k \geq k_0$.
Then, by our choice of $r$ in \eqref{eq:rchoice}, we have that
\begin{equation*}
  \frac{\vert \Sigma \vert}{2^{n-1} \, i(\Sigma) \, \vert \Sp^{n-1}\vert \, r_0^{n-1}} \leq k
  \iff r \leq \frac{r_0}{2}.
\end{equation*}
Hence for $k \geq k_0$, \eqref{eq:croke} holds.

We now define a test function $g_i$ supported on $B(x_i,2r)$ as follows
\begin{equation*}
  g_i(x) =
  \begin{cases}
    1, & \mbox{if } x \in B(x_i,r), \\
    1 - \frac{d(x,B(x_i,r))}{r}, & \mbox{if } x \in B(x_i,2r) \setminus B(x_i,r),\\
    0, & \mbox{otherwise}.
  \end{cases}
\end{equation*}
Now
\begin{equation*}
  \int_{\Sigma} u_i^2 \, d\V_{\Sigma} \geq \int_{\Sigma \cap B(x_i,r)} u_i^2 \, d\V_{\Sigma}
  = \vert \Sigma \cap B(x_i,r) \vert \geq \frac{2^{n-1} \vert \Sp^{n-1} \vert^n}{(n-1)^{n-1} \vert \Sp^n \vert^{n-1}} \, r^{n-1}.
\end{equation*}
In addition, $\vert \nabla g_i \vert^2 \leq \frac{1}{r^2}$ almost everywhere so
\begin{equation*}
  \int_{M} \vert \nabla g_i \vert^2 \, d\V_{M} = \int_{M \cap B(x_i,2r)} \vert \nabla g_i \vert^2 \, d\V_{M}
  \leq \frac{\vert M \cap B(x_i,2r) \vert}{r^2},
\end{equation*}
and by \eqref{eq:m} we obtain
\begin{equation*}
  \int_{M} \vert \nabla g_i \vert^2 \, d\V_{M} \leq 2^{n-1} i(M) \, \vert \Sp^n \vert \, r^{n-2}.
\end{equation*}
Hence the Rayleigh quotient $R(g_i)$ satisfies
\begin{equation*}
  R(g_i) \leq \frac{(n-1)^{n-1}  \vert \Sp^n \vert^n \, i(M)}{ \vert \Sp^{n-1} \vert^n} \cdot \frac{1}{r}.
\end{equation*}
We deduce that for $k \geq k_0$,
\begin{equation*}
  \sigma_k(M) \leq \tilde{B}(n) \, i(M) \, \left(\frac{i(\Sigma) k}{\vert \Sigma \vert}\right)^{1/(n-1)},
\end{equation*}
where $\tilde{B}(n) = \frac{4 (n-1)^{n-1}  \vert \Sp^n \vert^n}{ \vert \Sp^{n-1} \vert^{n - (1/(n-1))}}$.

Now by \eqref{eq:k0}, we have
\begin{equation*}
  k_0 \leq \frac{2 \vert \Sigma \vert}{2^{n-1} \, i(\Sigma) \, \vert \Sp^{n-1}\vert \, r_0^{n-1}}
\end{equation*}
Hence for $k \leq k_0$,
\begin{equation*}
  \sigma_k(M) \leq \sigma_{k_0}(M) \leq \tilde{A}(n) \frac{i(M)}{r_0},
\end{equation*}
where $\tilde{A}(n) = \frac{2^{1 + (1/(n-1))} (n-1)^{n-1}  \vert \Sp^n \vert^n}{ \vert \Sp^{n-1} \vert^{n}}$.

Hence for any $k \geq 1$, we conclude that
\begin{equation*}
  \sigma_k(M) \leq \tilde{A}(n)\frac{i(M)}{r_0} + \tilde{B}(n) \, i(M)\left(\frac{i(\Sigma) k}{\vert \Sigma \vert}\right)^{1/(n-1)}
\end{equation*}
as required.
\end{proof}

\begin{rem}
For the special case where $\Sigma \subset \R^n$ is an $(n-1)$-dimensional, closed, smooth hypersurface bounding a domain $M \subset \R^n$, we can use the inequality
\begin{equation*}
 \vert M \cap B(x,r)\vert \leq \omega_n r^{n}
\end{equation*}
instead of Inequality \eqref{eq:m} in the above proof. Then $i(M)$ is no longer present in the upper bound as stated in Corollary \ref{cor:hyp}.
\end{rem}

\section{Example to show that the injectivity radius is necessary}\label{s:4}

The goal of this section is to show that the term involving the injectivity radius in
Inequality \eqref{eq:T3} of Theorem \ref{thm:3} is necessary in the case where the
dimension of $M$ is at least 4. To do this, for $n \geq 3$,
we construct an $(n+1)$-dimensional submanifold of $\R^{n+3}$ whose boundary has volume $1$ and such
that the intersection index of the manifold and its boundary are bounded, but the first non-zero
Steklov eigenvalue of the manifold tends to $+\infty$ as the injectivity radius tends to zero.

The idea of our construction is as follows. We begin with an $n$-dimensional annulus in $\R^n$
which has a Steklov boundary condition on the interior (smaller) boundary sphere and a Neumann boundary
condition on the exterior (larger) boundary sphere. Roughly speaking, if the radius of the interior boundary sphere is small and the radius of the exterior boundary sphere is large, then the first non-zero eigenvalue of the mixed Steklov--Neumann problem on the annulus is large (see Subsection \ref{ss:4.1}).

In Subsection \ref{ss:4.2}, we construct a manifold $M$ and show that it has large first non-zero Steklov eigenvalue. The idea of the construction is that $M$ contains a \emph{nice} set for which it is possible to estimate the first non-zero Steklov eigenvalue of the mixed Steklov--Neumann problem. The latter gives a lower bound for the first non-zero Steklov eigenvalue of $M$ by the bracketing \eqref{eq:SN}. This estimate is quite technical and we defer it to Subsection \ref{calculation}.

More precisely, this \emph{nice} set is the Cartesian product of the annulus considered in Subsection \ref{ss:4.1} with a circle of a certain radius such that the volume of the interior boundary is equal to $1$. Allowing the width of the annulus to become large and the radius of the interior boundary sphere to become small gives rise to large
eigenvalues.

A difficulty is then to ``close up'' the annulus so that we can view the product
as a submanifold $M$ of $\R^{n+3}$ with one boundary component. This is done in Subsection \ref{construction}.
Note that it is also necessary to control the index of $M$. This is the object of Subsection \ref{index}.

\subsection{First non-zero eigenvalue of an annulus with mixed Steklov--Neumann boundary condition}\label{ss:4.1}~\\
We note that similar computations to those carried out in this subsection can be found in \cite{FS}. These calculations will be used at the end of Subsection \ref{ss:4.2} (see Inequality \eqref{eq:4.3}).

For $0 < \epsilon < \delta$, we consider the annulus
$$ A(\epsilon, \delta) = \{ x \in \R^n : \epsilon < d(x,0) < \delta\}.$$
Let $\partial A_{\delta}$ denote the boundary component of radius $\delta$ and $\partial A_{\epsilon}$
denote the boundary component of radius $\epsilon$. We have that $\vert \partial A_{\epsilon} \vert = n\omega_n\epsilon^{n-1}$. 
We consider the following eigenvalue problem on $A = A(\epsilon, \delta)$ which has a Neumann boundary condition on $\partial A_{\delta}$ and a Steklov boundary condition on $\partial A_{\epsilon}$:
\begin{gather*}
    \Delta u =0 \mbox{ in } A \\
    \partial_{\nu} u =0 \mbox{ on } \partial A_{\delta} \\
    \partial_{\nu} u = \sigma u \mbox{ on } \partial A_{\epsilon}.
\end{gather*}
\indent
We consider spherical coordinates $\rho$, $\theta_1, \dots, \theta_{n-1}$ where $\rho$ denotes the radial variable.
In these coordinates, it is well known that the Laplacian can be written as follows
\begin{equation*}
  \Delta = \frac{\partial^2}{\partial \rho^2} + \frac{(n-1)}{\rho}\frac{\partial}{\partial \rho}
  + \frac{1}{\rho^2} \Delta_{\Sp^{n-1}},
\end{equation*}
where $\Delta_{\Sp^{n-1}}$ denotes the Laplacian on the $(n-1)$-dimensional unit sphere $\Sp^{n-1}$.
Below we use that the corresponding eigenvalues of $\Sp^{n-1}$ are $k(k+n-2)$.
\newline
\indent
We now separate the angular and radial variables and consider
$$u(\rho,\theta) = \alpha(\rho) \beta(\theta_1, \dots, \theta_{n-1}).$$
We have that
\begin{equation*}
  \Delta u = \left(\alpha''(\rho) + \frac{(n-1)}{\rho}\alpha'(\rho) - \frac{\alpha(\rho)}{\rho^2} k(k+n-2)\right)
  \beta(\theta_1, \dots, \theta_{n-1}).
\end{equation*}
So if $\Delta u = 0$ then, as $\beta(\theta_1, \dots, \theta_{n-1})$ is not identically zero,
\begin{equation*}
  \left(\alpha''(\rho) + \frac{(n-1)}{\rho}\alpha'(\rho) - \frac{\alpha(\rho)}{\rho^2} k(k+n-2)\right)
  =0.
\end{equation*}
By setting $\alpha(\rho) = \rho^{\theta}$, we deduce that $\theta = k$ or $\theta = -n +2 -k$.
Hence
\begin{equation*}
  \alpha(\rho) = a\rho^{k} + b\rho^{-n + 2 - k}, \quad a,b \in \R,
\end{equation*}
with the boundary conditions
\begin{equation*}
  \alpha'(\delta) = 0, \quad \quad -\alpha'(\epsilon) = \sigma \alpha(\epsilon).
\end{equation*}
From $\alpha'(\delta) = 0$, we obtain that
\begin{equation*}
  b = \frac{a k \delta^{2k-2+n}}{k-2+n}.
\end{equation*}
By substituting this formula for $b$ into $-\alpha'(\epsilon) = \sigma \alpha(\epsilon)$, we obtain
\begin{equation*}
  -k\epsilon^{k-1} + k \delta^{2k-2+n}\epsilon^{-k+1-n}
  =\sigma\left(\epsilon^{k} + \frac{k}{k-2+n} \delta^{2k-2+n}\epsilon^{-k+2-n}\right).
\end{equation*}
In what follows, we will be interested in the first non-zero eigenvalue, so we set $k=1$ and obtain
\begin{equation*}
  -1 + \delta^{n}\epsilon^{-n}
  =\sigma\left(\epsilon + \frac{1}{n-1} \delta^{n}\epsilon^{1-n}\right),
\end{equation*}
which implies that
\begin{equation}\label{eq:4.1}
  \sigma_1^{N}(A(\epsilon,\delta)) = \frac{-1 + \delta^{n}\epsilon^{-n}}{\epsilon + \frac{1}{n-1} \delta^{n}\epsilon^{1-n}}.
\end{equation}

\subsection{Example to show that the injectivity radius is necessary}\label{ss:4.2}~\\
The idea is to construct a family of manifolds $M= M_{\epsilon,\delta,R}$ with boundary of volume $1$ such that the injectivity radius of the boundary tends to $0$ as $\epsilon$ tends to $0$.

We first consider the product $  A(\epsilon,\delta) \times S_R^1 \subset \R^n \times \R^2$, where $S_R^1$ is a circle of radius $R$, and view it as follows:
$$ A(\epsilon,\delta) \times S_R^1 \subset \R^n \times \R^2 \subset \R^{n+2} \times \{0\} \subset \R^{n+3}.$$
We note that $\delta$, $R$ depend on $\epsilon$.
We have $\Sp_{\epsilon}^{n-1} = \partial A_{\epsilon} \subset \R^{n-1}$.

We would like to view $A(\epsilon,\delta) \times S_R^1$ as part of an $(n+1)$-dimensional submanifold $M= M_{\epsilon,\delta,R}$ of $\R^{n+3}$ with boundary isometric to $\Sigma=\Sp^{n-1}_{\epsilon} \times S_R^1$.
The difficulty is to ``close up'' the annulus $A(\epsilon,\delta)$, that is to see it as part of a manifold $N=N_{\epsilon,\delta}$ with boundary $\Sp^{n-1}_{\epsilon}$. We also need to be able to show that the intersection index of the manifold that we construct (and that of its boundary) is bounded from above. In fact, we can construct a manifold of revolution in $\mathbb R^{n+1}$ with one boundary component.
\medskip

\subsubsection{Construction of the manifold}\label{construction}~\\
We consider $N=N_{\epsilon,\delta}$ as the union of three parts:
\begin{itemize}
\item the annulus $A(\epsilon,\delta) \subset \R^n$ that we view in $\R^{n+1}$ as
$$ \tilde{N}_1 := A(\epsilon, \delta) \times \{x_{n+1} = -1\}, $$
\item the ball $B_{\delta} \subset \R^n$ of centre $0$ and radius $\delta$ that we view in $\R^{n+1}$ as
$$ \tilde{N}_2 := B_{\delta} \times \{x_{n+1} = 1\},$$
\item the manifold of revolution given by the equation
$$ x_1^2 + \dots + x_n^2 = (\delta + \sqrt{1 - x_{n+1}^2})^2, \, -1 \leq x_{n+1} \leq 1,$$
which is contained in the real algebraic variety of degree $4$:
$$ \tilde{N}_3 := \{(x_1^2 + \dots + x_{n}^2 -\delta^2 -1 + x_{n+1}^2)^2 - 4 \delta^2(1-x_{n+1})^2 = 0,
\, -1 \leq x_{n+1} \leq 1\}.$$
\end{itemize}

The manifold $N$ is only $C^1$ but we can make it $C^2$ using a manifold of revolution given by an
equation of degree 4.
We then consider the submanifold
$$ M= M_{\epsilon,\delta,R} = N_{\epsilon,\delta} \times S_{R}^1 \subset \R^{n+1} \times \R^2$$
which is of dimension $n + 1$, has boundary $\Sp_{\epsilon}^{n-1} \times S_{R}^1$ and contains $A(\epsilon,\delta) \times S_{R}^1$.
We note that $\partial _{I} A(\epsilon,\delta)$ is smooth.
We see that $\vert \Sp_{\epsilon}^{n-1} \times S_R^1 \vert = 2\pi n\omega_n R \epsilon^{n-1}$.
To ensure that $\vert \Sp_{\epsilon}^{n-1} \times S_R^1 \vert = 1$, we choose $R = (2\pi n \omega_n)^{-1}\epsilon^{-(n-1)}$.

\medskip
\subsubsection{Intersection index of the manifold and its boundary}\label{index}~\\
We now show that the intersection index of $M =N \times S_R^1$, respectively of $\Sigma=\Sp_{\epsilon}^{n-1} \times S_R^1$, is bounded.
In particular, as $\Sigma = \Sp_{\epsilon}^{n-1} \times S_R^1 \subset \R^{n+3}$ is $n$-dimensional, we have to show that
$i_3(\Sp_{\epsilon}^{n-1} \times S_R^1)$ is bounded.
Similarly, as $N \times S_R^1 \subset \R^{n+3}$ is $(n+1)$-dimensional, we have to show that
$i_2(N \times S_R^1)$ is bounded.

We wish to apply Lemma 1 of \cite{M64} to obtain an upper bound for the intersection index of each of the
manifolds $\tilde{N}_1 \times S_{R}^1$, $\tilde{N}_2 \times S_{R}^1$ and $\tilde{N}_3 \times S_{R}^1$.
Lemma 1 of \cite{M64} states that it is possible to bound the intersection index of a real algebraic variety
defined by polynomial equations from above by the product of the degrees of these polynomials.
Then, as $M \subset \cup_{i=1}^{3} \tilde{N}_i \times S_{R}^1$, taking the sum of these
upper bounds gives an upper bound for the intersection index of $M$.

\smallskip
Let $Z_i = \tilde{N}_i \times S_{R}^1$ for $i=1,2,3$. Let $m \geq n+3$ and $p \geq 2$.
Let $X \subset \R^m$ be a $p$-plane. Let $S \subset \R^m$ be an open ball of $\R^m$.
Define $F: X \times S \to \R^m$ by $F(x,s) = x + s$.
Then by the Transversality Theorem (see, for example, \cite{GP74}), if $F$ is transverse to $Z_i$,
then for almost every $s \in S$, $f_s(x) = F(x,s)$ is transverse to $Z_i$.
That is, almost every $p$-plane is transverse to $Z_i$.
Hence by Lemma 1 of \cite{M64}, $i_p(Z_1) \leq 1 \times 2 = 2$, $i_p(Z_2) \leq 1 \times 2 = 2$
and $i_p(Z_3) \leq 4 \times 2 = 8$. So $i_p(Z_1 \cup Z_2 \cup Z_3) \leq 12$ and
$i_p(M) \leq 12$. As $\Sigma=\partial M = \Sp_{\epsilon}^{n-1} \times S_R^1 \subset Z_1 \times S_{R}^1$, $i_p(\partial M) \leq 2$.

\subsubsection{Steklov problem on the manifold}\label{calculation}~\\
We consider the Steklov problem on $M$:
\begin{gather*}
    \Delta u =0 \mbox{ in } N \times S_R^1 \\
    \partial_{\nu} u = \sigma u \mbox{ on } \Sigma=\partial (N \times S_R^1).
\end{gather*}
We are interested in the non-zero Steklov eigenvalues of $M$.
For $x \in N$, $y \in S_R^1$, we take a separation of variables
to obtain
$$ u(x,y) = f(x) g(y)$$
where $g$ is an eigenfunction of the Laplacian on $S_R^1$ with corresponding eigenvalue $\lambda \geq 0$.
Then
\begin{equation*}
  0 = \Delta u = (\Delta f) g + f(\Delta g) = (\Delta f - \lambda f)g.
\end{equation*}
In addition, for each fixed $y \in S_R^1$, we have $\nabla_{x} u (x,y) = (\nabla_{x} f(x)) g(y)$ and $\nu_{\Sigma }(x,y) = \nu(x)$, where $\nu(x) \in \R^n \times \{0\} $ denotes the outward unit normal to the boundary $\Sp_{\epsilon}^{n-1}$ of $N$. So, for each fixed $y \in S_R^1$,
\begin{equation*}
  \frac{\partial u}{\partial \nu_{\Sigma}}(x,y) = \nabla_{x} u(x,y) \cdot \nu_{\Sigma}(x,y)
  = (\nabla_{x} f(x) \cdot \nu(x)) g(y) =\left( \frac{\partial f}{\partial \nu}(x)\right) g(y).
\end{equation*}
Hence we have that
\begin{gather*}
    \Delta f = \lambda f \mbox{ in } N \\
    \frac{\partial f}{\partial \nu} = \sigma f \mbox{ on } \Sp_{\epsilon}^{n-1}.
\end{gather*}
We observe that if $\lambda = 0$, then $\sigma$ is a non-zero Steklov eigenvalue of $N$.
Indeed, if $\lambda = 0$, then the first eigenfunction $g$ is constant and non-zero.
As $u$ is a Steklov eigenfunction of $M$
with corresponding eigenvalue $\sigma > 0$, we have
\begin{equation*}
  0 = \int_{\Sigma} u \, d\V_{\Sp_{\epsilon}^{n-1}} \, d\V_{S_R^1} = \int_{\Sp_{\epsilon}^{n-1}} f \, d\V_{\Sp_{\epsilon}^{n-1}} \int_{S_R^1} g \, d\V_{S_R^1},
\end{equation*}
so $\int_{\Sp_{\epsilon}^{n-1}} f \, d\V_{\Sigma} =0$, i.e. $f$ is orthogonal to constants on $\Sp_{\epsilon}^{n-1}$.
\newline
\indent
On one hand, we have that
\begin{equation*}
  \int_{N \times S_R^1} f \Delta f \, d\V_{N} \, d\V_{S_R^1} = \V(S_R^1) \int_{N} f \Delta f \, d\V_{N}
  =\lambda \vert S_R^1 \vert \int_{N} f^2 \, d\V_{N}.
\end{equation*}
On the other hand, by Green's formula, we have that
\begin{align*}
  \int_{N \times S_R^1} f \Delta f \, d\V_{N} \, d\V_{S_R^1} &= \vert S_R^1 \vert \int_{N} f \Delta f \, d\V_{N} \\
  &= \vert S_R^1 \vert \left(\int_{\Sp_{\epsilon}^{n-1}} \partial_{\nu}f \, f \, d\V_{\Sp_{\epsilon}^{n-1}} - \int_{N} \vert \nabla f \vert^2 \, d\V_{N}\right)\\
  &= \vert S_R^1 \vert \left(\sigma \int_{\Sp_{\epsilon}^{n-1}} f^2 \, d\V_{\Sp_{\epsilon}^{n-1}} - \int_{N} \vert \nabla f \vert^2 \, d\V_{N}\right).
\end{align*}
So we obtain
\begin{equation*}
  \sigma = \frac{\lambda \int_{N} f^2 \, d\V_{N} + \int_{N} \vert \nabla f \vert^2 \, d\V_{N}}{\int_{\Sp_{\epsilon}^{n-1}} f^2 \, d\V_{\Sp_{\epsilon}^{n-1}}}.
\end{equation*}
Hence if $\lambda = 0$, then
\begin{equation}\label{eq:4.2}
  \sigma = \frac{\int_{N} \vert \nabla f \vert^2 \, d\V_{N}}{\int_{\Sp_{\epsilon}^{n-1}} f^2 \, d\V_{\Sp_{\epsilon}^{n-1}}}\geq \sigma_1(N) \geq \sigma_1^{N}(A(\epsilon,\delta)),
\end{equation}
where we use that $f$ is orthogonal to constants on $\Sp_{\epsilon}^{n-1}$ in the first inequality, and Steklov--Neumann bracketing \eqref{eq:SN} in the second inequality.
\\[10pt]
We first treat the case where $\lambda > 0$ and we come back to the case where $\lambda = 0$ at the end.
Let $A = A(\epsilon,\delta) \subset N$ where $\delta > \epsilon$ is to be chosen below. We have that
\begin{equation*}
  \sigma = \frac{\lambda \int_{N} f^2 \, d\V_{N} + \int_{N} \vert \nabla f \vert^2 \, d\V_{N}}{\int_{\Sp_{\epsilon}^{n-1}} f^2 \, d\V_{\Sp_{\epsilon}^{n-1}}} \geq \frac{\lambda \int_{A} f^2 \, d\V_{N} + \int_{A} \vert \nabla f \vert^2 \, d\V_{N}}{\int_{\Sp_{\epsilon}^{n-1}} f^2 \, d\V_{\Sp_{\epsilon}^{n-1}}}.
\end{equation*}
On $A$ we consider the metric
\begin{equation*}
  \tilde{g}(r,q) = dr^2 + h(r)^2 g_0
\end{equation*}
where $\epsilon \leq r \leq \delta$, $q \in \Sp^{n-1}$, $h(r) = r$ and $g_0$ is the canonical metric on $\Sp^{n-1}$.
Let $\{v_i\}_{i=0}^{\infty}$ be an orthonormal basis of eigenfunctions of the Laplacian on $\Sp^{n-1}$ with corresponding eigenvalues $\mu_i$. Any smooth function $f \in L^2(A)$ can be written in the following form:
\begin{equation*}
  f(r,q) = \sum_{i=0}^{\infty} a_i(r) v_i(q).
\end{equation*}
We have that
\begin{equation*}
  \int_{\Sp_{\epsilon}^{n-1}} f^2 \, d\V_{\Sp_{\epsilon}^{n-1}} = \Vert f \Vert_{L^2(\Sp_{\epsilon}^{n-1})}^{2}
  = \sum_{i=0}^{\infty} a_i(\epsilon)^2 h(\epsilon)^{n-1},
\end{equation*}
and
\begin{equation*}
  \int_{A} f^2 \, d\V_{N} = \Vert f \Vert_{L^2(A)}^{2}
  = \sum_{i=0}^{\infty} \int_{\epsilon}^{\delta} a_i(r)^2 h(r)^{n-1} \, dr.
\end{equation*}
In addition, $df(r,q) = \sum_{i=0}^{\infty} (a_i'(r) v_i(q) dr + a_i(r) dv_i(q))$ so
\begin{align*}
  \int_{A} \vert \nabla f\vert^2 \, d\V_{M} = \Vert \nabla f \Vert_{L^2(A)}^{2}
  &= \sum_{i=0}^{\infty} \int_{\epsilon}^{\delta} \left(a_i'(r)^2 + \frac{a_i(r)^2 \mu_i}{h(r)^2}\right) h(r)^{n-1} \, dr \\
  &= \sum_{i=0}^{\infty} \int_{\epsilon}^{\delta} \left(a_i'(r)^2 h(r)^{n-1} + a_i(r)^2 \mu_i h(r)^{n-3}\right)\, dr.
\end{align*}
For each $i \in \N$, we consider the ratio
\begin{equation*}
  R_i = \frac{\int_{\epsilon}^{\delta} \left(a_i'(r)^2 h(r)^{n-1} + \lambda a_i(r)^2 h(r)^{n-1} + a_i(r)^2 \mu_i h(r)^{n-3}\right)\, dr}{a_i(\epsilon)^2 h(\epsilon)^{n-1}}.
\end{equation*}
For each $\epsilon \leq r_1 \leq \delta$,
\begin{gather*}
  \int_{\epsilon}^{r_1} a_i'(r) \, dr = a_i(r_1) - a_i(\epsilon), \\
  \bigg\vert \int_{\epsilon}^{r_1} a_i'(r) \, dr \bigg\vert^2 \leq (r_1 - \epsilon) \int_{\epsilon}^{r_1} a_i'(r)^2 \, dr.
\end{gather*}
So
\begin{equation*}
  \int_{\epsilon}^{r_1} a_i'(r)^2 \, dr \geq \frac{(a_i(r_1) - a_i(\epsilon))^2}{r_1 - \epsilon},
\end{equation*}
and, as $r \mapsto h(r)$ is an increasing function,
\begin{equation*}
  \int_{\epsilon}^{r_1} a_i'(r)^2 h(r)^{n-1} \, dr \geq h(\epsilon)^{n-1}\frac{(a_i(r_1) - a_i(\epsilon))^2}{r_1 - \epsilon}.
\end{equation*}
First suppose there exists $\epsilon \leq r_1 \leq \delta$ such that $\vert a_i(r_1)\vert < \frac{\vert a_i(\epsilon)\vert}{2}$. Then
\begin{equation*}
  \int_{\epsilon}^{\delta} a_i'(r)^2 h(r)^{n-1} \, dr \geq
  \int_{\epsilon}^{r_1} a_i'(r)^2 h(r)^{n-1} \, dr
  \geq h(\epsilon)^{n-1}\frac{a_i(\epsilon)^2}{4(r_1 - \epsilon)}.
\end{equation*}
Hence we obtain that
\begin{equation*}
  R_i \geq \frac{\int_{\epsilon}^{\delta} a_i'(r)^2 h(r)^{n-1} \, dr}{a_i(\epsilon)^2 h(\epsilon)^{n-1}} \geq \frac{1}{4(r_1 - \epsilon)}.
\end{equation*}
We choose $r_1 = 2\epsilon$ then $r_1 - \epsilon >0$ and $R_i \geq \frac{1}{4 \epsilon}$.
\newline
On the other hand, if $\vert a_i(r)\vert \geq \frac{\vert a_i(\epsilon)\vert}{2}$ for all $\epsilon \leq r \leq \delta$, then
\begin{align*}
  R_i &\geq \frac{\int_{\epsilon}^{\delta} (\lambda a_i(r)^2 h(r)^{n-1} + a_i(r)^2 \mu_i h(r)^{n-3}) \, dr}{a_i(\epsilon)^2 h(\epsilon)^{n-1}}
  \geq \frac{\int_{\epsilon}^{\delta} (\lambda h(r)^{n-1} + \mu_i h(r)^{n-3}) \, dr}{4 h(\epsilon)^{n-1}}.
\end{align*}
If $i \geq 1$, then $\mu_i \geq n-1$ and
\begin{equation*}
  R_i \geq \frac{\int_{\epsilon}^{r_1} (n-1) r^{n-3} \, dr}{4 \epsilon^{n-1}}
  = \frac{(n-1)}{4(n-2)} \left(\frac{(2\epsilon)^{n-2} - \epsilon^{n-2}}{\epsilon^{n-1}}\right)
  = \frac{(n-1)}{4(n-2)} \left(\frac{2^{n-2} - 1}{\epsilon}\right).
\end{equation*}
If $i=0$, then $\mu_0 = 0$ and
\begin{align*}
  R_i &\geq \frac{\int_{\epsilon}^{\delta} \lambda a_i(r)^2 h(r)^{n-1} \, dr}{a_i(\epsilon)^2 h(\epsilon)^{n-1}}
  \geq \frac{\int_{\epsilon}^{\delta} \lambda r^{n-1} \, dr}{4 \epsilon^{n-1}}
  = \frac{\lambda}{4n} \left(\frac{\delta^n - \epsilon^n}{\epsilon^{n-1}}\right).
\end{align*}
As $\lambda \geq \lambda_1(S_R^1) = \frac{1}{4R^2} = \pi^2 n^2 \omega_n^2 \epsilon^{2(n-1)}$,
\begin{equation*}
  R_i \geq  \frac{n \pi^2 \omega_n^2}{4} \epsilon^{2(n-1)} \left(\frac{\delta^n - \epsilon^n}{\epsilon^{n-1}}\right).
\end{equation*}
We choose $\delta = 2\epsilon^{-1}$. Then $r_1 = 2\epsilon < \delta$ as $\epsilon <1$, and we have
\begin{equation*}
  R_i \geq \frac{n \pi^2 \omega_n^2}{4} \left(\frac{2^n - \epsilon^{2n}}{\epsilon}\right)
  \geq \frac{n \pi^2 \omega_n^2 (2^n -1)}{4} \cdot \frac{1}{\epsilon}.
\end{equation*}
Hence, we have that
\begin{equation*}
  R_i \geq \frac{\tilde{C}}{\epsilon}
\end{equation*}
for $i \in \N$ where
\begin{equation*}
  \tilde{C} = \min \left\{ \frac{1}{4}, \frac{(2^{n-2}-1)(n-1)}{4(n-2)}, \frac{n \pi^2 \omega_n^2 (2^n -1)}{4}\right\},
\end{equation*}
so $\sigma \geq \frac{\tilde{C}}{\epsilon}$.
\\[10pt]
\indent
We now use the choice of $\delta = 2\epsilon^{-1}$ to deal with the case where $\lambda=0$.
By \eqref{eq:4.2} and \eqref{eq:4.1} from Subsection \ref{ss:4.1}, we have that
\begin{equation}\label{eq:4.3}
  \sigma \geq \sigma_1^{N}(A(\epsilon,2\epsilon^{-1}))
  = \frac{-\epsilon^{2n} + 2^n}{\epsilon^{2n+1} + \frac{2^n}{n-1}\epsilon}
  \geq \frac{(2^n -1)}{(1+\frac{2^n}{n-1})\epsilon}
  =\frac{(n-1)(2^n -1)}{(n-1 +2^n)}\cdot \frac{1}{\epsilon}.
\end{equation}

We conclude that
\begin{equation*}
  \sigma \geq \frac{\hat{C}}{\epsilon}
\end{equation*}
where
\begin{equation*}
  \hat{C} = \min \left\{ \frac{1}{4}, \frac{(2^{n-2}-1)(n-1)}{4(n-2)}, \frac{n \pi^2 \omega_n^2 (2^n -1)}{4}, \frac{(n-1)(2^n -1)}{(n-1 +2^n)} \right\}.
\end{equation*}
Therefore $\sigma \to +\infty$ as $\epsilon \to 0$.

\section{An obstruction to the optimal exponent}\label{s:5}
In this section we present two examples of Riemannian manifolds $(M,g)$ for which the presence of
$\vert M \vert_{g}^{\beta}$, $\beta >0$, in an upper bound for $\sigma_k(M,g)$ poses an obstruction to achieving the optimal exponent of $k$.

Let $(M,g)$ be a compact Riemannian manifold of dimension $n \geq 2$ with smooth boundary $\Sigma$.
We often have the following type of inequality for the Steklov eigenvalues of $M$: for each $k\geq 1$,
$$
\sigma_k(M,g)\leq C(n) \vert \Sigma\vert^{\gamma}\vert M\vert^{\beta}\, k^{\alpha},
$$
where $\alpha >0$ and $\beta \geq 0$. In order for this inequality to be invariant by homothety, we require that $\gamma=-\frac{1+n\beta}{n-1}$.

A particular inequality that we wish to obtain is the case where $\alpha=\frac{1}{n-1}$ and $\beta=0$.
However, a simple example shows that $\alpha$ and $\beta$ cannot be chosen independently.
In fact, for every such manifold $M$, we can construct an example of a metric $g$ such that $\alpha$ and $\beta$ satisfy the inequality
\begin{equation}\label{eq:2.1}
 1+\beta \leq \alpha (n-1).
\end{equation}

For example, if $\beta=0$, then $\alpha \geq \frac{1}{n-1}$, and if $\beta>0$, then $\alpha > \frac{1}{n-1}$.
In the latter case, the bound cannot have the optimal exponent of $k$ with respect to the Steklov spectral asymptotics \eqref{eq:1.1}. In general, the only situation where we can hope to obtain a bound with optimal exponent of $k$ is when $\beta =0$, that is, there is no contribution from the volume of $M$.

We first consider the particular case of a cylinder and then address a more general situation.
In what follows, we make use of Lemma 2.1 and Inequality $(5)$ from \cite{CGG}.

\begin{ex}\label{ex:2.1}
Let $\Sigma$ be a smooth, closed, connected Riemannian manifold. Let $L > 0$ and $M = \Sigma \times [0,L]$.
If, for $k \geq 1$, we have
$$
\sigma_k(M,g)\leq C(\vert \Sigma\vert,n) \vert M\vert^{\beta}\, k^{\alpha},
$$
where $\alpha >0$, $\beta \geq 0$, and $C(\vert \Sigma\vert,n)$ is a constant depending only on $\vert \Sigma\vert$, $n$, then \eqref{eq:2.1} holds.

The Steklov spectrum of $M = \Sigma \times [0,L]$ is the union of $0$, $2/L$, $\sqrt {\lambda_k}\tanh(\sqrt {\lambda_k}L/2)$ and $\sqrt {\lambda_k}\coth(\sqrt {\lambda_k}L/2)$, where $\lambda_k$ are the Laplace eigenvalues of $\Sigma$.
In particular, $\sigma_{2k}(M) \geq \sqrt {\lambda_k}\tanh(\frac{\sqrt {\lambda_k}L}{2})$.
Then taking $L=1/\sqrt{\lambda_k}$, this implies
$$
\sqrt{\lambda_k} \tanh \frac{1}{2} \le \sigma_{2k}(M) \le C(\vert \Sigma \vert,n) \vert M\vert^{\beta}(2k)^{\alpha},
$$
and we have
$$
\lambda_k^{\frac{1+\beta}{2}}\le \tilde C(\vert \Sigma\vert,n) k^{\alpha}.
$$
As $\lambda_k \sim k^{2/(n-1)}$ as $k \to \infty$, we conclude that inequality \eqref{eq:2.1} holds.
\end{ex}

\begin{ex}\label{ex:2.2}
For a compact Riemannian manifold $(M,g)$ of dimension $n \geq 2$ with smooth boundary $\Sigma$, we can construct an example of a metric $g$ such that if
$$
\sigma_k(M,g)\leq C(\vert \Sigma\vert_{g},n) \vert M\vert_{g}^{\beta}\, k^{\alpha},
$$
where $\alpha >0$, $\beta \geq 0$, and $C(\vert \Sigma\vert_{g},n)$ is a constant depending only on $\vert \Sigma\vert_{g}$, $n$, then \eqref{eq:2.1} holds.
\newline
\indent
Suppose $M$ has $b$ boundary components denoted $\Sigma_1,\dots,\Sigma_b$. We can construct a Riemannian metric $g$ on $M$ such that near each boundary component $\Sigma_i$, $i=1,\dots,b$, this metric is isometric to a product $\Sigma_i \times [0,L_i]$.
\newline
\indent
Indeed, by compactness, there exists $L>0$ such that around each $\Sigma_i$, $M$ is diffeomorphic to $\Sigma_i \times [0,2L]$ via a diffeomorphism $\phi_i$, and the $\Omega_i=\phi_i^{-1}(\Sigma_i \times [0,2L])$ are disjoint.
On each neighbourhood $\Omega_i$, there are two different metrics: the restriction of $g$ and the metric $g_i$ such that $(\Omega_i,g_i)$ is isometric to the product metric on $\Sigma_i \times [0,2L]$.
We introduce a smooth decreasing function $\psi$ on $[0,2L]$ that takes the value $1$ on $[0,L]$ and the value $0$ around $2L$. We then consider the new metric
	$$
	g_{\psi,i}= \psi g_i+\epsilon^2(1- \psi)g
	$$	
This is the product metric on $\Sigma_i \times [0,L]$ and the metric $\epsilon^2 g$ near $2L$. We can hence extend it by $\epsilon^2 g$ outside $\cup_{i=1}^b \Omega_i$, so the volume of the complement of $\cup_{i=1}^b \Omega_i$ in $M$ is as small as we like. We note that throughout the process described above, the metric does not change near the $\Sigma_i$ hence the volume of each boundary component $\Sigma_i$ does not change either.
\newline
\indent
By Lemma 2.2 of \cite{CGG}, the mixed Steklov-Neumann spectrum of a cylinder $\Sigma \times [0,L]$ of length $L$ is given by $\sigma_k^N(\Sigma \times [0,L])=\sqrt {\lambda_k}\tanh(\frac{\sqrt {\lambda_k}L}{2})$  where $\lambda_k$ are the Laplace eigenvalues of $\Sigma$.
In particular, $\sigma^N_{k}(\Sigma \times [0,L]) \ge \sqrt {\lambda_k}\tanh(\frac{\sqrt {\lambda_k}L}{2})$.
\newline
\indent
We set $L_i=\frac{1}{\sqrt{\lambda_k(\Sigma_i)}}$, and suppose that $k$ is large enough such that $L_i\leq L$ for $i=1,\dots,b$.
The choice of $L_i$ implies that
$$
\sqrt{\lambda_k(\Sigma_i)} \tanh \frac{1}{2} \le \sigma^N_{k} (\Sigma_i\times [0,L_i]).
$$
\indent
Let $\bar{\sigma_k}(M,g)$ be the minimum of $\sigma_{[\frac{k}{b}]}(\Sigma_i\times [0,L_i])$ where $[.]$ denotes the integer part. Then, we have
$$
\sigma_k(M,g)\ge \bar{\sigma_k}(M,g).
$$
\indent
So, for each given $k$, there exists $i(k)\in 1,\dots,b$ such that
$$
\sigma_k(M,g) 
\ge \sigma_{[\frac{k}{b}]}(\Sigma_{i(k)}\times [0,L_{i(k)}]) \ge
\sigma_{[\frac{k}{b}]}^{N}(\Sigma_{i(k)}\times [0,L_{i(k)}]) \ge
\tanh \frac{1}{2}\sqrt{\lambda_{[\frac{k}{b}]}(\Sigma_{i(k)})},
$$
where, in the second inequality, we use Steklov-Neumann bracketing \eqref{eq:SN}.
But, for each $i=1,\dots,b$, $\lambda_{[\frac{k}{b}]}(\Sigma_i) \sim [\frac{k}{b}]^{2/(n-1)}\sim k^{2/(n-1)}$ as $k \to \infty$.
\newline
\indent
We note that by the choice of the metric $g_{\psi,i}$, only the volume of the cylindrical neighbourhood of $\Sigma$ plays a role in what follows.
\newline
\indent
On one hand, we have
\begin{equation}\label{eq:2.2}
\tanh \frac{1}{2}\sqrt{\lambda_{[\frac{k}{b}]}(\Sigma_{i(k)})} \le \sigma_k(M,g) \le  C(\vert \Sigma\vert,n) \vert M\vert_{g_{\psi,i}}^{\beta}k^{\alpha} \le C(\vert \Sigma\vert,n) k^{\alpha} \vert \Sigma \vert^{\beta}(\sum_{i=1}^b L_i)^{\beta}.
\end{equation}
\indent
On the other hand, we have
$$
(\sum_{i=1}^b L_i)^{\beta}= \left(\sum_{i=1}^b \frac{1}{\sqrt {\lambda_{[\frac{k}{b}]}(\Sigma_i)}}\right)^{\beta}.
$$
\indent
For each $i$, as $k\to \infty$, we have $\lambda_{[\frac{k}{b}]}(\Sigma_i)\sim k^{2/(n-1)}$. (At this stage, it is crucial that the boundary components $\Sigma_1,\dots,\Sigma_b$ are fixed throughout the process). This implies that
$$
(\sum_{i=1}^b L_i)^{\beta}\sim k^{-\beta/(n-1)},
$$
which together with \eqref{eq:2.2} implies, as for the cylinder, that for large $k$
$$
k^{1/(n-1)} k^{\beta/(n-1)} \le C k^{\alpha},
$$
and we conclude that \eqref{eq:2.1} holds.
\end{ex}

\bigskip\noindent{\bf Acknowledgements.}
The authors wish to thank Luigi Provenzano for very useful remarks.

\bibliographystyle{plain}
\bibliography{bibliocolboisgittins}

\end{document}